\documentclass[s5,11pt]{article}
\usepackage{amsmath}
\usepackage{latexsym}
\newtheorem{theorem}{Theorem}

  \author{ Tord Sj\"odin, Ume\aa}
\title{ A Short and Unified Proof of Kummer's Test\footnote{An extended version of the paper was presented at the "Workshop on Beauty and Explanation in Mathematics",  Ume\aa\ University, Sweden, 2014} }
 \begin{document}  
\maketitle
\begin{abstract}Kummer's test from 1835 states that the positive series $\sum _{n=1}^\infty a_n$ is convergent if and only if there is a sequence $\{ B_n\}_1^\infty$ of positive numbers such that 
 $B_n\cdot  \frac{a_n }{a_{n+1}} -B_{n+1}\geq 1 ,$ for all  sufficiently large $n$. We present an exact analysis and a short and unified proof of Kummer's test. The test has been applied to differential equations and studied in mathematical philosophy.\end{abstract}
\paragraph{ \it AMS 2010 Subject Classsification:} Primary 40A05,
Secondary 00A30, 01A55
\paragraph{ \it Key words and phrases:}   Sequence, positive series, convergence, convergence test, Kummer's test, mathematical explanation \section{Introduction} 
Kummer's test appeared in the famous {\it Crelles Journal} in 1835, \cite{Ku} and is reprinted in his Collected Works, \cite{We} p. 47. Recall that a positive series is convergent if and only if the sequence of partial sums has a finite upper bound.  
\begin{theorem}A positive series $\sum\limits _1^\infty a_n$ is convergent if and only if there is a sequence $\{ B_n\}_1^\infty$ of positive numbers and an integer $N\geq 1$ such that 
 \begin{equation}B_n\cdot  \frac{a_n }{a_{n+1}} -B_{n+1}\geq 1 ,\quad \textrm{for all } n\geq N.\end{equation}
 \end{theorem}
Kummer's original proof was complicated and not readily accepted by his contemporaries. Fifty years later, Stoltz gave a proof of todays standard, see Knopp \cite{Kn}  p. 311, which is the proof seen today. It treats the necessity part and the sufficiency part in separate cases. We present a short  and unified proof of Kummer's test (Theorem 1) in the next section, together with some comments on how the test has been applied in mathematics and mathematical philosophy.
\section{Proof of Kummer's test}
   Our proof is very direct. We first note that it is no loss of generality to assume that (1) holds with equality and find an explicit formula for the numbers $\{ B_n\}_N^\infty$. This will allow us to give a very  short and unified proof.     \\[0.5em]
{\it Proof of Kummer's test.}   
We begin by rewriting (1) as
\begin{equation}B_{n+1}\leq B_n\cdot  \frac{a_n }{a_{n+1}}-1 ,\quad n\geq N. \end{equation}
Since all the numbers $B_n$ are positive, nothing is lost by increasing $B_{n+1}$ at each step in (2) such that we get an equality
\begin{equation}B_{n+1}= B_n\cdot  \frac{a_n }{a_{n+1}}-1 ,\quad n\geq N.\end{equation}
 Then conditions (1)$-$(3) are equivalent in the sense that if there exists a positive sequence $\{ B_n\}_N^\infty$ satisfying anyone of them, then the same is true for the other two.
 Now we leave the positivity of $B_n$ aside and solve the recursion formula (3) in the standard way. Multiplying (3) by $a_{n+1}$ gives
  $$a_{n+1}\cdot B_{n+1}=a_n\cdot B_n- a_{n+1},\quad n\geq N.$$ 
  Adding these equations by telescoping and rearranging the terms we get the formula
  \begin{equation}
 B_{n+1}=\frac{1}{a_{n+1}}\cdot \bigg(  B_N\cdot a_N -  (a_{N+1}+\cdots +a_{n+1}) \bigg), \quad n\geq N,\end{equation}
 which is our basis for the proof.  \\[0.5em]
 First assume that the series is convergent with sum $\sum_{n=1}^\infty a_n=s $. Choose $N=1$ and any $B_1>s/a_1$, then the sequence $\{ B_n\}_1 ^\infty $ constructed by (4) is positive and satisfies (3) and (1). Conversely, assume that there is a positive sequence $\{ B_n\}_N^\infty $ satisfying (1). Then the sequence $\{B_n\}_N^\infty $ constructed from (3) and (4) is positive as well and the series is convergent since the number $B_N\cdot a_N$ is an upper bound for the partial sums $a_{N+1}+a_{N+2}+\cdots +a_{n+1}$, for all $n\geq N$.\hfill $\Box$\\[0.5em]
 Our proof of Kummer's test shows that the sequence $\{ B_n\}_1^\infty$ in (1) is almost uniquely determined by the terms of the series itself, in contrast to all earlier proofs.
Kummer's test is not of much help in proving convergence for positive series in general, since there is no rule how to find $N$ and $B_N$ in (4). It can however be quite useful in some special cases. If we rewrite (1)  as
\begin{equation}\frac{a_{n+1}}{a_ n}\leq \frac{B_n}{\rho +B_{n+1}}, \quad \textrm{for all}\,  n\geq N \textrm{ and some $\rho >0$ }. \end{equation}
it becomes a rate of decrease for the terms in the series comparable to some of the classical convergence tests, such as the Quotient Test, Root Test and Gauss Test, \cite{D}. See also \cite{GH}. where (5) is applied to difference equations and \cite{HM}, \cite{Ste}, where the test is used in the study of mathematical explanations. 
\\[0.5em]
  {\it Remark.} Eduard Kummer was born in Prussia (todays Germany) in 1810 and was appointed to a chair in Berlin in 1855. He is best remembered for his ideal numbers and ideals, that paved the way for modern ring theory. We think that our proof reveals the true nature of the test and that the sequences $\{ B_n\}_N^\infty$ constructed by (3) and satisfying (4) might be of practical use.\ 

\end{document}